\documentclass[letterpaper, 10 pt, conference]{ieeeconf}  

\IEEEoverridecommandlockouts                              

\usepackage{amsfonts}
\usepackage{amsmath}
\usepackage{algorithm}
\usepackage{algorithmicx,algpseudocode}
\usepackage{graphicx}
\usepackage{epstopdf}

\usepackage{amsthm}
\usepackage{booktabs}

\newtheorem{assump}{Assumption}
\newtheorem{lem}{Lemma}
\newtheorem{thm}{Theorem}
\newtheorem{rem}{Remark}

\DeclareMathOperator{\bfz}{\mathbf{z}}

\title{\LARGE \bf
An Improved Primal-Dual Interior Point Solver\\ for Model Predictive Control*
}

\author{Xi Zhang$^{1}$, Laura Ferranti$^{2}$ and Tam\'as Keviczky$^{2}$
\thanks{*This research is supported by the European Union's Seventh Framework
	Programme FP7/2007-2013 under grant agreement n. AAT-2012-RTD-2314544
	(RECONFIGURE) and by the TU Delft Space Institute.}
\thanks{$^{1}$
        {\tt\small cindyzhang0919@gmail.com}}%
\thanks{$^{2}$Delft University of Technology, Delft, 2628 CD, The Netherlands
        {\tt\small \{l.ferranti,t.keviczky\}@tudelft.nl}}%
}
\newif\ifpaper 
\paperfalse

\begin{document}

\maketitle
\thispagestyle{empty}
\pagestyle{empty}

\begin{abstract}
We propose a primal-dual interior-point (PDIP) method for solving
quadratic programming problems with linear inequality constraints that typically arise
form MPC applications. We show that the solver converges (locally) quadratically
to a suboptimal solution of the MPC problem. PDIP solvers rely on two phases:
the damped and the pure Newton phases. Compared to state-of-the-art PDIP
methods, our solver replaces the initial damped Newton phase (usually
used to compute a medium-accuracy solution) with a dual solver based on Nesterov's fast gradient scheme (DFG) 
that converges with a sublinear convergence rate of order
$\mathcal{O}(1/k^2)$ to a medium-accuracy solution.
The switching strategy to the pure Newton phase, compared to the state of the art, is computed in the dual space to exploit the
dual information provided by the DFG in the first phase.
Removing the damped Newton phase has the additional advantage that our solver
saves the computational effort required by backtracking line search. The effectiveness of the proposed solver is demonstrated on a 2-dimensional discrete-time unstable system and on an aerospace application.

\end{abstract}


\section{Introduction}

Model predictive control (MPC) is an advanced control technique that offers an
elegant framework to solve a wide range of control problems
(regulation, tracking, supervision, etc.) and handle constraints on the plant.
The control objectives and the constraints are usually formulated as an
optimization problem that the MPC controller has to solve (either offline or
online) to return the control command for the plant. In this work, we focus on
MPC problems with quadratic cost and linear constraints that can be formulated
as quadratic programming (QP) problems (which cover a large range of
practical applications).

The presence of this optimization problem has traditionally limited the use of
MPC to \emph{slow} processes, that is, processes with no hard real-time constraints.
Recently, MPC has received increasing attention in fields, such as
aerospace and automotive, where the real-time aspects are
critical and the computation time for the controller is limited. 
Hence, offline (e.g., the explicit MPC~proposed
by~\cite{bemporad2002explicit}) and online (using solvers tailored for MPC
applications) solutions have been investigated to overcome the computational
issues related to the MPC controller. In this work, we focus on online
solutions that allow one to handle a wider range of problems. Online
optimization algorithms can be divided in two main families: first- and second-order methods.
First-order methods, such as gradient or splitting methods and their
accelerated
versions~\cite{nesterov1983method,kogel2011fast,nesterov2013introductory,Stathopoulos2016}, have simpler 
theoretical requirements (such as, Lipschitz continuity only on
the first derivative of the cost) and converge to a
medium-accuracy solution within few iterations. Second-order methods, such as 
active-set \cite{ricker1985ASMuse,schmid1994ASMquadratic,ferreau2006ASMonline} 
and interior-point methods 
\cite{rao1998application,wang2010fast,zeilinger2014realtimeRobustMPC}, have
more strict theoretical requirements (such as, Lipschitz continuity on
the first and second derivative of the cost), but are more suitable when a
high-accuracy solution is required. In this work, we are mainly interested in
primal-dual interior-point (PDIP) solvers~\cite{Borrelli2015predictive}. As
detailed for example in~\cite{Borrelli2015predictive,boyd2004convex}, their
convergence can be divided in two phases: (i) \emph{damped Newton phase}
characterized by a linear convergence rate and used to reach a
medium-accuracy solution; (ii) \emph{pure Newton phase} characterized by a
quadratic convergence rate and used to improve the accuracy of the solution
obtained from phase (i).

\emph{\textbf{Contribution}.} The main contribution of this paper is a
PDIP solver for solving
inequality-constrained QP problems that commonly arise
from MPC applications. The proposed solver combines the advantages of Nesterov's fast gradient (FG) 
method~\cite{nesterov2013introductory} and of PDIP solvers. In particular, we exploit the ability of the fast 
gradient method to converge to medium-accuracy solutions within few
iterations (that can be performed efficiently) and the ability of PDIP solvers to converge to high-accuracy 
solutions in Phase (ii). The proposed combination aims to improve the convergence of the PDIP by replacing Phase 
(i) with the FG method that has a sublinear convergence rate of order $\mathcal
O(1/k^2)$.
We modify the classical analysis of the PDIP~\cite{boyd2004convex} to take into account the presence of (active) inequality constraints and allow the switch between the two solvers. This is done by moving the convergence analysis of the solver and the formulation of the switching strategy to the dual framework. As a consequence, we provide bounds on the level of primal suboptimality and infeasibility achieved with our proposed algorithm. An additional feature of the proposed solver is that the computational effort related to the backtracking line search (required in Phase (i), refer to \cite[Chapter 9]{boyd2004convex} in details) is removed, given that Phase (ii) uses a unit step size.
Finally, the proposed approach is tested on two MPC applications, that
is, the regulation of an unstable input- and output- constrained planar system and the Cessna Citation Aircraft system \cite{maciejowski2002predictive}.

\emph{\textbf{Notation}.} We work in the space $\mathbb{R}^n$ composed of column
vectors. $\left\Vert \bfz\right\Vert_1$ and $\left\Vert
\bfz\right\Vert_2$
indicate the 1-norm and 2-norm, respectively.
$[\bfz]_+$ and $[\bfz]_-$ indicate the projection onto nonnegative
orthant and negative orthant, respectively. Furthermore, $\operatorname{diag}(\bfz)$ indicates a matrix that has the elements of $\bfz$ on its main diagonal. Finally, $\mathbf 1$ is the vector of ones.  
\section{Problem Formulation}
\label{sec:problem_formulation}

In this work we focus on the control of \emph{discrete, linear, time-invariant} systems, which can be described as follows:
\begin{subequations}
\label{subeq:2_plant}
	\begin{align}
		x({t+1})& = Ax(t)+Bu(t)\\
		y(t) & = Cx(t)+Du(t)
	\end{align}
\end{subequations}
where $x\in\mathcal{X}\subseteq\mathbb{R}^{n_x}$, $u\in\mathcal{U}\subseteq\mathbb{R}^{n_u}$, $y\in\mathcal{Y}\subseteq \mathbb R^{n_y}$ denote the state, control command, output, respectively. The sets $\mathcal{U}$, $\mathcal{X}$, $\mathcal{Y}$ are closed convex sets that contain the origin in their interior. In addition, $t\geq 0$ denotes the sampling instant. Furthermore, the system matrices $A$, $B$, $C$, and $D$ are constant matrices of fixed dimension. In the remainder of the paper, we assume that the pair $(A,B)$ is stabilizable.

In this work, we focus on regulation problems that can be addressed by solving online 
(i.e., every time new measurements $x(t)$ are available from the plant) the following optimization problem:
\begin{subequations}
\label{eq:2_QP_MPC}
	\begin{align}
		J^*= \underset{x_k\in \mathcal{X},
		u_k\in\mathcal{U}}{\operatorname{minimize}}~~ &
		\sum_{k=0}^{N}J_k(x_k,u_k) \\
		\textrm{subject to}:~&x_{k+1}= Ax_{k}+Bu_{k}\label{eq:plant_dyn_prediction}\\
		& y_k \in \mathcal Y,~k = 0,...,N,\\
		&~ x_0 = x(t), 
	\end{align}
\end{subequations}
where $x_k$, $u_k$, and $y_k$ represent the predicted evolution of the state, 
control command, and output, respectively. Furthermore, $J_k(x_k,u_k)$ denotes
the stage cost and is defined as follows:
	\begin{equation*}
		J_k(x_k,u_k) = \frac 1 2\begin{cases} \left(x_k^{\textrm{T}}Qx_k + u_k^{\textrm{T}}Ru_k\right) & k =
		0,...,N-1,\\
		x_N^{\textrm{T}}Qx_N & k=N,
	\end{cases}
	\end{equation*}
where $Q=Q^{\textrm{T}}\succeq 0$ and $R=R^{\textrm{T}}\succ 0$ weigh the state and input, respectively. 

\textbf{\emph{Condensed Formulation}.} We exploit the plant dynamics~\eqref{eq:plant_dyn_prediction}
 to eliminate the states from the decision variables and express them as an explicit function of the current measured
  state $x(t)$ and future control inputs 
  (refer to \cite{maciejowski2002predictive} for more details). This is known as \emph{condensed} QP formulation, which leads to
   compact and dense QPs with the control inputs as decision variables:
\begin{equation}
	\mathbf z:=\begin{bmatrix}
		u_0^{\textrm{T}}&u_{1}^{\textrm{T}}&\ldots&u_{N-1}^{\textrm{T}}
	\end{bmatrix}^{\textrm{T}}.
	\label{eq:2_condensed}
\end{equation}
This leads to a set of equalities with dense matrices expressing the decision variables as a function of the current state and input sequence:
\begin{subequations}
	\begin{equation}
		\mathbf{x}=A_N{x(t)}+B_N\mathbf{z}
		\label{subeq:2_condensed_state}
	\end{equation}
	\begin{equation}
		\mathbf {y}=C_N{x(t)}+D_N\mathbf{z}
		\label{subeq:2_condensed_output}
	\end{equation}
\end{subequations}
where $\textbf{x}=[x_0^{\textrm{T}}~x_1^{\textrm{T}}~\ldots~x_{N}^{\textrm{T}}]$ and $\textbf{y}=[y_0^{\textrm{T}}~y_1^{\textrm{T}}~\ldots~y_{N}^{\textrm{T}}]$. For details on the structure of $A_n$, $B_n$, $C_n$ and $D_n$ refer to \cite{maciejowski2002predictive}.

Further, define $\mathcal{Z}:=\mathcal U\times \ldots\times \mathcal U$ ($N$ times).
Finally, the optimization problem in \eqref{eq:2_QP_MPC} can be written in the following standard form which only consists of inequality constraints:
\begin{subequations}
\label{eq:2_QP_standard}
	\begin{align}
		\min_{\bfz}~&~f_0(\bfz):=\frac{1}{2}\bfz^{\textrm{T}}H\bfz+(h{x(t)})^{\textrm{T}}\bfz\\
		s.t.~&~g(\bfz):=G\bfz+E{x(t)}+g\le 0\label{eq:ineq_constr_condensed}
	\end{align}
\end{subequations}
where $H\succ0$ (given $R\succ0$), $h$, $G$, $E$, and $g$ are given matrices 
(more details on their structure can be found
in~\cite{maciejowski2002predictive}) and \eqref{eq:ineq_constr_condensed} encodes $\mathcal U$, $\mathcal X$, and $\mathcal Y$.

In this paper, we consider the following assumptions:
\begin{assump}
\label{ass:ass_strong_convexity}
	Function $f_0(\bfz)$ is $m_p$-strongly convex and twice differentiable, $m_p$ is the convexity parameter.
\end{assump}
\begin{assump}
\label{ass:slater}
	The Slater condition holds for Problem \eqref{eq:2_QP_standard}, i.e., there exists ${\bfz}\in\operatorname{relint}(\mathcal{Z})$ with $g(\bfz)<0$.
\end{assump}
\section{Preliminaries}
\label{sec:preliminaries}
\subsection{Dual Fast Gradient Method}
\label{subsec:dfg}
In the following, we provide a high level description of the fast gradient method proposed by Nesterov~(refer to~\cite{nesterov1983method,nesterov2005smooth,nesterov2013introductory} for more details) with a focus on its application to the dual of Problem~\eqref{eq:2_QP_standard}, as Algorithm~\ref{alg:dual_fast_gradient} details. 
In particular, notice that Algorithm~\ref{alg:dual_fast_gradient} exploits (steps 5 and 6) a projection step at each iteration. 
If this projection is hard to compute, this operation can be challenging to
solve efficiently online. This is the case, if we apply
Algorithm~\ref{alg:dual_fast_gradient} directly to Problem~\eqref{eq:2_QP_standard}, in which the set $g(\bfz)\le0$ consists of so-called complicating constraints. Hence, we operate on the dual of Problem~\eqref{eq:2_QP_standard} that only requires the computation of the projection on the positive orthant, which can be computed efficiently. The dual of Problem~\eqref{eq:2_QP_standard} can be computed as follows. Given the Lagrangian described below:
\begin{equation}
\mathcal{L}(\bfz,\lambda)=f_0(\bfz)+\lambda^{\textrm{T}}g(\bfz)
\label{eq:3_lagrangian}
\end{equation}
where $\lambda\in\mathbb{R}^m$ is the Lagrange multiplier, the dual function is described as follows:
\begin{equation}
d(\lambda)=\min_{\bfz}\mathcal{L}(\bfz,\lambda)
\label{eq:3_dual_function}
\end{equation}
which first-order derivative is given as follows 
\[\nabla d(\lambda)=g(\bfz).\]
The dual of Problem~\eqref{eq:2_QP_standard} is given by 
\begin{equation}
f_0(\bfz^*)=d(\lambda^*)=\max_{\lambda\geq 0}d(\lambda),
\label{eq:3_outer_problem}
\end{equation}
where the $f_0(\bfz^*)\!=\!d(\lambda^*)$ follows from
Assumptions~\ref{ass:ass_strong_convexity}~and~\ref{ass:slater}.
\begin{algorithm}[t]
\fontsize{8.5}{8.5}\selectfont
\caption{Dual Fast Gradient Method.}
	\label{alg:dual_fast_gradient}
	\begin{algorithmic}[1] 
		\State{Given $H$, $h$, $g$, $E$, $G$, $x(t)$, $\hat \lambda$, $k_{\max}$, $L_d$.}
		\State Initialize $\lambda_0=\hat \lambda$. 
		\For{$k=0,\ldots, k_{\max}$}
		\State Compute $\bfz_{k}=\arg\min_{\bfz}\mathcal{L}(\bfz,\lambda_k)$.\label{alg:dfg_1}
		\State Compute $\hat{\lambda}_k=\left[\lambda_k+\frac{1}{L_d}\nabla d(\lambda_k)\right]_+$.\label{alg:dfg_2}
		\State Compute $\lambda_{k+1}=\frac{k+1}{k+3}\hat{\lambda}_k+\frac{2}{L_d(k+3)}\left[\sum_{j=0}^{k}\frac{j+1}{2}\nabla d(\lambda_j)\right]_+.$\label{alg:dfg_3}
		\EndFor
	\end{algorithmic}
\end{algorithm} 

Algorithm~\ref{alg:dual_fast_gradient} uses $1/L_d$ ($L_d:=\left\Vert
GH^{-1}G^{\textrm{T}}\right\Vert_2$) as step size. T
his step size is optimal for the proposed algorithm~\cite{giselsson2014improved}. In particular, as 
detailed in \cite{giselsson2014improved},  $d(\lambda)$ is Lipschitz continuous with constant $L_d$, leading 
to a tighter upper bound than the one provided for example in \cite{nesterov2013introductory,necoara2014rate}. 
Furthermore, the algorithm initializes $\lambda_0 = \hat \lambda\ge 0$ (step 2). 
A complete analysis for Algorithm \ref{alg:dual_fast_gradient} starting from $\hat \lambda\neq 0$ can be found in \cite{necoara2014rate}. For simplicity, we start Algorithm \ref{alg:dual_fast_gradient} with $\hat \lambda=0$. At every iteration, the algorithm first computes a minimizer for Problem~\eqref{eq:3_dual_function} (step 4). Then, it performs a linear update of the dual variables (steps 5-6).
 
\textbf{\emph{Convergence Analysis}.} Algorithm~\ref{alg:dual_fast_gradient} has
a convergence rate of $\mathcal{O}(1/k^2)$, as demonstrated in~\cite{nesterov1983method,nesterov2013introductory}. In addition, 
as discussed for example in~\cite{necoara2014rate}, Algorithm \ref{alg:dual_fast_gradient} converges to a suboptimal solution of MPC problem after $k_{\max}$ iterations (refer to~\cite{necoara2014rate} for more details).
In the remainder of the paper, we make the following assumption:
\begin{assump}
\label{ass:5_final_active_set}
	There exists a $k_{\max}\geq 0$ such that after $k_{\max}$ iterations, Algorithm~\ref{alg:dual_fast_gradient} is able to find a solution $\bfz_{k_{\max}}$ close to the central path, such that, for $\eta_d\geq 0$, the following holds:
	\begin{equation*}
	\left\Vert\left[g(\bfz_{k_{\max}})\right]_+\right\Vert_2\le\eta_d.
	\end{equation*}
\end{assump} 
\subsection{Primal-Dual Interior Point Method}
\label{subsec:pdip}
\begin{algorithm}[t]
\fontsize{8.5}{8.5}\selectfont
	\caption{Primal-Dual Interior Point Method.}
	\label{alg:primal_dual_interior_point}
	\begin{algorithmic}[1] 
	    \State{Given $H$, $h$, $g$, $E$, $G$, $x(t)$.}
		\State Initialize $\bfz_0$, $\lambda_0>0$, $s_0>0$, centering parameter $\kappa\in(0,1)$, backtracking line search parameters $\alpha\in(0,0.5),~\beta\in(0,1)$, tolerance $\varepsilon>0$, and $k=0$.
		\Repeat
		\State Determine $\tau_{k+1}=\mu_{k+1}=\kappa\mu_k$.\label{alg:pdip_1}
		\State Compute search direction $\Delta \zeta_{\textrm{pd}}$ by solving \eqref{eq:3_search_direction_infeasible}.\label{alg:pidp_2}
		\State Backtracking line search $\rho_k:=1$
		\While {$f_0(\bfz +\rho_k\Delta\bfz)>f_0(\bfz)+\alpha\rho_k\nabla f_0(\bfz)^T\Delta\bfz$}
		\State $\rho_k:=\beta\rho_k$.
		\EndWhile
		\State Update $\zeta_{k+1}=\zeta_k+\rho_k\Delta \zeta_{\textrm{pd}}$ where $\rho_k>0$ is the step size.
		\State $k = k+1.$ 
		\Until {stopping criterion $\mu_k\le\varepsilon$.}\label{alg:pdip_3}
		\State \textbf{return} Point close to $\bfz^*$ from $\zeta_k=(\bfz_k,\lambda_k,s_k)$.
	\end{algorithmic}
\end{algorithm}
In the following, we present a version of the primal-dual interior point (PDIP) method proposed by~\cite{Borrelli2015predictive} and described in Algorithm~\ref{alg:primal_dual_interior_point}.
 
The general idea of primal-dual interior point methods is to solve the KKT conditions by using a modified version of Newton's method. In this respect, recall the Lagrangian defined in \eqref{eq:3_lagrangian}. PDIP solves the following \emph{relaxed} KKT conditions:
\begin{subequations}
\label{eq:3_rKKT}
	\begin{equation}
	\nabla f_0(\bfz)+Dg(\bfz)^{\textrm{T}}\lambda=0
	\label{subeq:3_rKKT_stationary}
	\end{equation}
	\begin{equation}
	g(\bfz)+s=0
	\label{subeq:3_rKKT_primal}
	\end{equation}
	\begin{equation}
	S\lambda=\tau\mathbf 1_m,
	\label{subeq:3_rKKT_centering}
	\end{equation}
	\begin{equation}
	(s,\lambda)>0,
	\label{subeq:3_rKKT_dual}
	\end{equation}
\end{subequations}
where $S:=\operatorname{diag}(s)$, $s\in\mathbb{R}^m$ is the slackness variable, $Dg(\bfz)$ 
is the derivative matrix of the inequality constraint function $g(\bfz)$, and $\tau \in[0,\mu]$, 
where $\mu=s^{\textrm{T}}\lambda/m$ denotes the average duality gap.

From \eqref{eq:3_rKKT}, the residual variable is defined as follows:
\begin{equation}
r_{\tau}(\bfz,\lambda,s)=\begin{bmatrix}
\nabla f_0(\bfz)+Dg(\bfz)^{\textrm{T}}\lambda\\
g(\bfz)+s\\
S\lambda-\tau\textbf{1}
\end{bmatrix}=\begin{bmatrix}
r_{\textrm{dual}}\\r_{\textrm{pri}}\\r_{\textrm{cent}}
\end{bmatrix}.
\label{eq:3_residual_infeasible}
\end{equation} 
The search direction can be obtained by
linearizing~\eqref{eq:3_residual_infeasible} at the current iterate
$\zeta_k=(\bfz_k,\lambda_k,s_k)$:
\begin{equation}
\underbrace{\begin{bmatrix}
	\nabla^2f_0(\bfz_k) & Dg(\bfz_k)^{\textrm{T}} & 0\\
	Dg(\bfz_k) & 0 & I\\
	0 & S_k & \Lambda_k
	\end{bmatrix}}_{Dr_{\tau}(\zeta)}\underbrace{\begin{bmatrix}
	\Delta \bfz_{\textrm{pd}}\\ \Delta\lambda_{\textrm{pd}}\\ \Delta s_{\textrm{pd}}
	\end{bmatrix}}_{\Delta\zeta_{\textrm{pd}}}=-\begin{bmatrix}
r_{\textrm{dual}}\\r_{\textrm{pri}}\\r_{\textrm{cent}}
\end{bmatrix}
\label{eq:3_search_direction_infeasible}
\end{equation}
where $\Lambda_k=\textbf{diag}(\lambda_k)$ and we use $\lambda^{\textrm{T}}\nabla^2g(\bfz)=0$.

\textbf{\emph{Convergence Analysis}.} The convergence of PDIP has been shown in \cite[Chapter 10]{boyd2004convex}. The convergence analysis is based on the convergence of the residual variable $r_{\tau}$ and can be divided into two phases. In the first phase, that is, the \emph{damped Newton phase}, Algorithm~\ref{alg:primal_dual_interior_point} converges \emph{linearly}, while in the second phase, that is, the \emph{pure Newton phase}, the backtracking line search selects unit step size and \emph{quadratic} convergence rate can be achieved. In order to enter pure Newton phase, the 2-norm of the residual variable has to satisfy the following condition
\begin{equation}
\left\Vert r_{\tau}(\bfz,\lambda,s)\right\Vert\le\eta_p
\label{eq:3_2nd_phase_condition}
\end{equation}
where $0<\eta_p\le\frac{{m}_{Dr}^2}{L_{Dr}}$ with ${m}_{Dr}$ the lower bound on $Dr_{\tau}(\zeta)$ and $L_{Dr}$ denotes the Lipschitz constant of $Dr_{\tau}(\zeta)$.

Loosely speaking, the condition above states that if the algorithm is far from the optimal solution it
converges more slowly (damped Newton phase), while when it is close to the optimal solution,
it converges faster (pure Newton phase). We rely on this observation to design an improved PDIP method that 
fully replaces the (slow) damped Newton phase with the dual fast gradient (DFG).
In particular, we propose to exploit the DFG
until~\eqref{eq:3_2nd_phase_condition} is satisfied, that is, we completely
remove the damped Newton phase in Algorithm~\ref{alg:primal_dual_interior_point}. Then, we initialize the 
PDIP with the solution returned by DFG to enter directly into the pure Newton phase, in which PDIP converges quadratically.
 As shown in the next section, combining these two methods is not trivial. One of the main issues is related to the presence of
  inequality constraints, compared to~\cite[Chapter 10]{boyd2004convex}, which only takes into account equality constraints. 
  In the presence of inequality,~\eqref{eq:3_2nd_phase_condition} (in the primal space) becomes too conservative to be used in practice (especially in the presence of active constraints), because the condition numbers of $Dr_{\tau}(\zeta)$ are hard to derive due to asymmetry and variance of $Dr_{\tau}(\zeta)$. Therefore, to overcome this issue, we propose, in the next section, a new switching condition (in the dual space) to enter the pure Newton phase.

\section{Proposed Solver}
\label{sec:proposed_solver}

As described in Section~\ref{subsec:pdip}, we would like to eliminate the damped Newton phase from PDIP and 
preserve the pure Newton phase, which allows the algorithm to converge quadratically to the optimal solution of Problem~\eqref{eq:2_QP_standard}.
Compared to the state of the art that analyzes the convergence of the algorithm
in the presence of equality constraints and in primal space, we move our analysis to the dual framework. This choice is strengthened
by the decision of using the DFG to replace the damped Newton phase. 
In particular, building on the convergence analysis in \cite[Chaper 9]{boyd2004convex} of the standard Newton's method,
we derive an estimate on primal suboptimality and feasibility violation achieved
with the proposed solver.

The following lemma (from~\cite[Theorem 2.1]{necoara2014rate}) is useful to prove the convergence of our proposed solver.

\begin{lem}
	Under Assumption \ref{ass:ass_strong_convexity} and \ref{ass:slater}, the dual function in \eqref{eq:3_dual_function} is twice differentiable. Then the gradient is given by:
	$\nabla d(\lambda)=g(\bfz(\lambda))$ and the dual Hessian is given by:
	\begin{equation}
	\label{eq:4_dual_hessian}
	\begin{aligned}
	\nabla^2d(\lambda)&=-\nabla g(\bfz(\lambda))\left[\nabla^2f_0(\bfz)\right]^{-1}\nabla g(\bfz(\lambda))^{\textrm{T}}\\
	&=-GH^{-1}G^{\textrm{T}}
	\end{aligned}
	\end{equation}
	Furthermore, the 2-norm of dual Hessian is bounded as follows:
	\begin{equation}
	{m}_d=\frac{\|G\|_2^2}{\sigma_{\max}(H)}\le\|\nabla^2d(\lambda)\|_2\le\frac{\|G\|_2^2}{\sigma_{\min}(H)}=M_d
	\label{eq:4_dual_H_bound}
	\end{equation}
	where $\sigma_{\max}(H)$ and $\sigma_{\min}(H)$ denote the maximum and minimum eigenvalues of $H\succ0$, respectively.
\end{lem}
\ifpaper
\begin{proof}
The proof can be found in~\cite{CDCArxiv}.
\end{proof}
\else
\begin{proof}
According to Assumption \ref{ass:ass_strong_convexity}, $f_0(\bfz)$ is $m_p-$strongly convex, that is, $H$ is positive definite, then the following holds:
\begin{equation}
\begin{aligned}
	\|\nabla^2d(\lambda)\|_2 &= \|GH^{-1}G^{\textmd{T}}\|_2\\
	& \ge\frac{1}{\sigma_{\text{max}}(H)}\|GIG^{\textmd{T}}\|_2
\end{aligned}
\end{equation}
using the fact that $H^{-1}=U\Sigma^{-1}U^{\textmd{T}}$ where $\Sigma$ is the singular value of $H$ and $U$ is an unitary matrix from singular value decomposition (SVD).
The proof of the upper bound is similar to the one in \cite[Theorem 2.1]{necoara2014rate}. In particular, we consider the derivative of the first-order optimality condition of Problem \eqref{eq:2_QP_standard} with respect to $\lambda$ to obtain $\nabla \bfz(\lambda)$. Then,~\eqref{eq:4_dual_hessian} follows by using $\lambda^{\textrm{T}}\nabla^2g(\bfz)=0$, given that only linear inequality constraints are involved in \eqref{eq:2_QP_standard}.
\end{proof}
\fi
From \eqref{eq:4_dual_hessian}, it follows that $\nabla^2 d(\lambda)$ is Lipschitz continuous with any $L_{dH}\ge0$.
\begin{thm}
	Under Assumptions~\ref{ass:ass_strong_convexity}-\ref{ass:5_final_active_set}, there exist $0<\eta_d\le{{m}_d^2}/{L_{dH}}$, $\gamma_d>0$ and $L_d>0$ such that the following holds:
	\begin{itemize}
		\item \emph{DFG Phase}. If
		$\left\Vert\left[g(\bfz)\right]_+\right\Vert_2\geq\eta_d$, DFG runs until the $\eta_d$-solution is achieved (after $k_{\max}$ iterations).
		\item \emph{Pure Newton Phase}. If
		$\left\Vert\left[g(\bfz)\right]_+\right\Vert_2\le\eta_d$, for all $k\geq 0$, with step size
		$\rho_k=1$, the following holds:
		\begin{equation}
		\left\Vert\nabla
		d(\lambda_{k+1})\right\Vert_2\le\frac{L_{dH}}{2{m}_d^2}\left\Vert\nabla
		d(\lambda_k)\right\Vert_2^2.
		\end{equation}
	     Hence, the following holds for $\bar k>0$:
		\begin{equation}
		0\geq f(\bfz_{\bar k})- f_0(\bfz^*)\geq -\left(\frac{{m}_d}{4L_{dH}}\right)^{\bar k}.
		\end{equation}
	\end{itemize}
	\label{thm:5_dual_convergence}
\end{thm}
\ifpaper
\begin{proof}
The proof can be found in~\cite{CDCArxiv}.
\end{proof}
\else
\begin{proof}
We summarize the proof here into the following steps:

\begin{itemize}
\item \emph{\textbf{Newton Increment}.} Similar to Newton decrement in primal space, in dual space, we define Newton increment for dual problem as following:
\begin{equation}
	\begin{aligned}
	\nu(\lambda)&=\left[-\nabla d(\lambda)^{\textmd{T}}\nabla^2d(\lambda)^{-1}\nabla d(\lambda)\right]^{\frac{1}{2}}\\
	& = \left[-\Delta\lambda_{\textmd{nt}}^{\textmd{T}}\nabla^2d(\lambda)\Delta\lambda_{\textmd{nt}}\right]^{\frac{1}{2}}
	\end{aligned}
\end{equation}
where $\Delta\lambda_{\textmd{nt}}=-\nabla^2d(\lambda)^{-1}\nabla d(\lambda)\ge0$ is the Newton search direction in dual space.

\item \emph{\textbf{Unit Step Size}.} In this step, we show that the following condition leads to a unit step size in backtracking line search in dual space:
\begin{equation}
	\left\Vert\nabla d(\lambda_k)\right\Vert_2\le\eta_d.
	\label{eq:4_2nd_phase_condition_dual}
\end{equation}
Note that due to space limitations, we omit some mathematical computations, which are similar to the convergence analysis proposed in \cite[Chapter 9]{boyd2004convex} for the proofs in primal space.
	
	By Lipschitz continuity of dual Hessian $\nabla^2d(\lambda)$, for $\rho\ge0$, the following relation holds:
	\begin{equation}
	\left\Vert\nabla^2d(\lambda+\rho\Delta\lambda_{\textmd{nt}})-\nabla^2d(\lambda)\right\Vert_2\le\rho L_{dH}\left\Vert\Delta\lambda_{\textmd{nt}}\right\Vert_2
	\label{eq:A_lipschitz_dual}
	\end{equation}
	
	Multiply \eqref{eq:A_lipschitz_dual} by $\Delta\lambda_{\textmd{nt}}$ on both sides and obtain:
	\begin{equation}
	\begin{aligned}
		\|\Delta\lambda_{\textmd{nt}}^T\left[\nabla^2d(\lambda+\rho\Delta\lambda_{\textmd{nt}})-\nabla^2d(\lambda)\right]\Delta\lambda_{\textmd{nt}}\|_2&\\
		\le\rho L_{dH}\left\Vert\Delta\lambda_{\textmd{nt}}\right\Vert_2^3&
	\end{aligned}
	\label{eq:A_lipschitz_dual+dlam}
	\end{equation}
	
	Define $\tilde{d}(\rho)=d(\lambda+\rho\Delta\lambda_{\textmd{nt}})$. Then $\nabla^2_{\rho}\tilde{d}(\rho)=\Delta\lambda_{\textmd{nt}}^T\nabla^2d(\lambda+\rho\Delta\lambda_{\textmd{nt}})\Delta\lambda_{\textmd{nt}}$. Therefore, \eqref{eq:A_lipschitz_dual+dlam} can be written as:
	\begin{equation}
	\left|\nabla^2_{\rho}\tilde{d}(\rho)-\nabla^2_{\rho}\tilde{d}(0)\right|\le \rho L_{dH}\left\Vert\Delta\lambda_{\textmd{nt}}\right\Vert_2^3
	\label{eq:A_lipschitz_dual_xi}
	\end{equation}
	Considering \eqref{eq:4_2nd_phase_condition_dual} and the strong convexity of function $-d(\lambda)$, by setting backtracking line search parameter $\alpha=\frac{1}{3}$ and integrating \eqref{eq:A_lipschitz_dual_xi} three times, we can prove the following:
	\begin{equation}
	d(\lambda+\Delta\lambda_{\textmd{nt}})\ge d(\lambda)+\alpha\nabla^2d(\lambda)^T\Delta\lambda_{\textmd{nt}}
	\label{eq:A_backtracking_result}
	\end{equation}
	which shows that with the unit step $\rho=1$ is accepted by backtracking line search in dual space.
	
\item \emph{\textbf{Modified Condition}.} Under Assumption \ref{ass:5_final_active_set}, \eqref{eq:3_dual_function} can be modified as follows:
\begin{equation}
	d(\lambda)\!=\!\min_{\bfz}
	\!f_0(\bfz)\!+\!\lambda^{\textrm{T}}g(\bfz)=\min_{\bfz}\!
	f_0(\bfz)\!+\!\lambda^{\textrm{T}}\left[g(\bfz)\right]_+
	\label{eq:5_dual_modified}
\end{equation}
where we use the fact that in $\lambda$ the entries corresponding to inactive constraints are zero. Therefore, the first-order derivative of the dual function $d(\lambda)$ for Problem \eqref{eq:5_dual_modified} is:
\begin{equation}
	\nabla d(\lambda)=\left[g(\bfz)\right]_+
\end{equation}
Thus, the switching condition~\eqref{eq:4_2nd_phase_condition_dual} becomes:
\begin{equation}
	\left\Vert\nabla d(\lambda_k)\right\Vert_2=\left\Vert\left[g(\bfz)\right]_+\right\Vert_2\le\eta_d
	\label{eq:4_condition_modified}
\end{equation}

\item \emph{\textbf{Upper Bound on Duality Gap}.} With the backtracking line search selecting unit step, by Lipschitz continuity, we can derive the following relation of the first-order derivative of the dual function between two iterations:
{\small
\begin{equation}
	\begin{aligned}
	\left\Vert\nabla d(\lambda_+)\right\Vert_2&=\left\Vert\int_{0}^{1}\left[\nabla^2d(\lambda+\rho\Delta\lambda_{\textmd{nt}})-\nabla^2d(\lambda)\right]\Delta\lambda_{\textmd{nt}}d\rho\right\Vert_2\\
	&\le\frac{L_{dH}}{2}\left\Vert\Delta\lambda_{\textmd{nt}}\right\Vert_2^2 = \frac{L_{dH}}{2}\left\Vert\nabla^2d(\lambda)^{-1}\nabla d(\lambda)\right\Vert_2^2\\
	& \le\frac{L_{dH}}{2m_d^2}\left\Vert\nabla d(\lambda)\right\Vert_2^2
	\end{aligned}
	\label{eq:A_gradient_2iter}
\end{equation}
}	
Therefore, for $\bar{k}>0$, recursively we can have
\begin{equation}
	\frac{L_{dH}}{2m_d^2}\left\Vert\nabla d(\lambda_{\bar{k}})\right\Vert_2\le\left(\frac{L_{dH}}{2m_d^2}\left\Vert\nabla d(\lambda)\right\Vert_2\right)^2
	\label{eq:A_quad_gradient}
\end{equation}
By relying on strong convexity property of function $-d(\lambda)$ (which follows from the assumptions of the theorem), a crucial result on the duality gap can be obtained, that is,
\begin{equation}
	f_0(\bfz^*) - f(\bfz_{\bar k})\leq
	f_0(\bfz^*)-d(\lambda_{\bar k})\le\frac{1}{2{m}_d}\left\Vert\nabla d(\lambda_{\bar k})\right\Vert_2^2
	\label{eq:4_bound_duality_gap}
\end{equation}
\end{itemize}

Finally, substituting~\eqref{eq:4_condition_modified}
into~\eqref{eq:4_2nd_phase_condition_dual}, \eqref{eq:4_2nd_phase_condition_dual} into \eqref{eq:A_quad_gradient}
and~\eqref{eq:A_quad_gradient} into \eqref{eq:4_bound_duality_gap}
concludes the proof.
\end{proof}
\fi
\begin{rem}
	The theorem above states that by relying on the DFG we can initialize the PDIP
	$\eta_d-$close to the central path. This allows the PDIP to "skip" the damped Newton phase
	and enter directly to pure Newton phase.
\end{rem}
\begin{rem}
Algorithm~\ref{alg:proposed_solver} does not
require backtracking line search given that the damped Newton phase has been
replaced by a less computational demanding DFG phase. In this respect,
note that the most computationally demanding step of the DFG is the (inner)
minimization of the Lagrangian, which can be computed efficiently and up to
a given accuracy (for details refer, for example, to
to~\cite{necoara2015adaptive} and the references within).
\end{rem}

Algorithm~\ref{alg:proposed_solver} summarizes the proposed strategy to compute a solution for the MPC problem~\eqref{eq:2_QP_MPC}, in which the DFG (steps 3-9) is used to compute an $\eta_d$-solution that allows the PDIP (steps 10-17) to enter directly in the pure Newton phase and converge quadratically to the optimal solution of the MPC problem.
When Algorithm~\ref{alg:proposed_solver} switches from the DFG phase to
the pure Newton phase (step 5), we have to make sure that we preserve the
information already computed by the DFG. Hence, the initialization strategy for
the PDIP is important to ensure a successful switch.  In this respect, step 5 of
Algorihtm~\ref{alg:proposed_solver} relies on the solution
$(\bfz_{\textrm{DFG}},\lambda_{\textrm{DFG}})$ returned by DFG as follows:
\begin{subequations}
\label{eq:4_initialization_strategy}
	\begin{align}
	\zeta_0&=(\bfz_{\textrm{DFG}},\,\lambda_{\textrm{DFG}},s_{0,{\textrm{PDIP}}}),\\
	s_{0,{\textrm{PDIP}}}&=-\left[g(\bfz_{\textrm{DFG}})\right]_-+\left[g(\bfz_{\textrm{DFG}})\right]_+.
	\label{subeq:initial_s}
	\end{align}
\end{subequations}
Equation~\eqref{subeq:initial_s} provides an initialization for $s$ that guarantees
$s>0$. Note that the information on the
duality gap is fully fed into the pure Newton phase by \eqref{subeq:initial_s},
according to the following
\begin{equation*}
	\begin{aligned}
	(s_{0,{\textrm{PDIP}}})^{\textrm{T}}\lambda_{0,\textrm{PDIP}}&=\left[-\left[g(\bfz_{\textrm{DFG}})\right]_-\!\!+\left[g(\bfz_{\textrm{DFG}})\right]_+\right]^{\textrm{T}}\!\!\lambda_{0,\textrm{PDIP}}\\
	&=\left[g(\bfz_{\textrm{DFG}})\right]_+^{\textrm{T}}\lambda_{0,\textrm{PDIP}},
	\end{aligned}
\end{equation*}
where we use the fact that $\left[g(\bfz_{\textrm{DFG}})\right]_-^{\textrm{T}}\lambda_{0,\textrm{PDIP}}=0$. 
\begin{algorithm}[t]
\fontsize{8.5}{8.5}\selectfont
	\caption{Proposed Solver.}
	\label{alg:proposed_solver}
	\begin{algorithmic}[1] 
	\State{Given $H$, $h$, $g$, $E$, $G$, $x(t)$, $\hat \lambda$, $L_d$, $\eta_d$, $\varepsilon<\eta_d$.}
		\State Initialize $\lambda_0=\hat \lambda$, $k=0$. 
		\Repeat
		\State Compute $\bfz_{k}=\arg\min_{\bfz}\mathcal{L}(\bfz,\lambda_k)$.
		\State Compute $\hat{\lambda}_k=\left[\lambda_k+\frac{1}{L_d}\nabla d(\lambda_k)\right]_+$.
		\label{alg:PS_step_hat_l}
		\State Compute $
		\lambda_{k+1}=\frac{k+1}{k+3}\hat{\lambda}_k+\frac{2}{L_d(k+3)}\left[\sum_{j=0}^{k}\frac{j+1}{2}\nabla d(\lambda_j)\right]_+.$
		\label{alg:PS_step_l}
		\State $k=k+1$.
		\Until {\begin{equation} \label{eq:4_condition_modified} \left\Vert\left[g(\bfz)\right]_+\right\Vert_2\le\eta_d
		\end{equation}}
		\State \textbf{return} $(\bfz_{\textrm{DFG}},\lambda_{\textrm{DFG}})$
		\State Initialize $\zeta_0$ according to \eqref{eq:4_initialization_strategy} and $k=0$.\label{eq:new_solver_2}
		\Repeat
		\State Determine $\tau_{k+1}=\mu_{k+1}=\kappa\mu_k=(s_k)^{\textrm{T}}\lambda_k/m$.
		\State Compute search direction $\Delta \zeta_{\textrm{pd}}$ by solving \eqref{eq:3_search_direction_infeasible}.
		\State Update $\zeta_{k+1}=\zeta_k+\Delta \zeta_{\textrm{pd}}$.
		\State $k=k+1$.
		\Until {stopping criterion $\mu_k\le\varepsilon$.}
		\State \textbf{return} Point close to $\textbf{z}^*$
	\end{algorithmic}
\end{algorithm}

\begin{rem}
	One could also try to derive the condition for PDIP to enter pure Newton phase in primal space directly when solving Problem 
	\eqref{eq:2_QP_standard}. In this respect, in primal space, the convergence of PDIP can be analyzed using its similarity
	 to the barrier method~\cite[Chapter 11]{boyd2004convex}. The main limitation, however, is that the algorithm has to be initialized 
	 with $r_{\textrm{pri}}=0$, that is, a strictly feasible starting point is
	 required.
	 This difficulty can be handled by either softening the constraints \cite{kerrigan2000soft} in the problem solved by the PDIP or by tightening the constraints \cite{necoara2015adaptive} in the problem solved by the DFG. These solutions are, however, more conservative. Furthermore, the switching condition in primal space involves the inverse of $g(\bfz)$ that becomes undefined when there are constraints active at the optimum.
\end{rem}

\section{Numerical Results}
\label{sec:evaluation}
The proposed solver (Algorithm \ref{alg:proposed_solver}) is tested both on the
planar discrete-time linear unstable system in
\cite{rubagotti2014stabilizing} and the open-loop unstable Cessna Citation Aircraft system in \cite{maciejowski2002predictive} (the descriptions of the two systems are omitted here for space limitations).

We tested the algorithm on a Windows OS, using an Intel(R) Core(TM) i7-4550 CPU (1.50-2.10 GHz) and RAM 8.00GB. The algorithms are implemented in Matlab.

\subsection{Planar Linear System}
The input and output of the system are subjected to the following constraints:
$\Vert u\Vert_{\infty}\le1,~\Vert y\Vert_{\infty}\le1$. Furthermore, $Q$ and $R$
are defined according to \cite{rubagotti2014stabilizing}.
As mentioned in Section \ref{sec:proposed_solver}, $L_{dH}$ can be selected as any
real number greater than zero. In this respect, we chose $L_{dH} = 200$.
Furthermore, ${m}_d =
1.4529$ and, consequently $\eta_d = 0.0106$. Finally, the backtracking line
search parameters for Algorithm~\ref{alg:primal_dual_interior_point} are chosen as
$\alpha=\frac{1}{3}$, $\beta=0.5$.

We compared the following 4 scenarios:
\begin{enumerate}
  \item Algorithm~\ref{alg:primal_dual_interior_point} is
  \emph{warm-started} with $\bfz_0:={\bfz}_{\textrm{LS}}$, which is the optimal
  solution of the unconstrained problem, that is, a least-squares (LS) problem,
  associated with Problem \eqref{eq:2_QP_standard}, and $\lambda_0 = \mathbf 1_m$.
\item Algorithm~\ref{alg:primal_dual_interior_point} is \emph{warm-started} with $\bfz_0:={\bfz}_{\textrm{LS}}$ as in Scenario 1 and $\lambda_0 =
10^{-6}\mathbf 1_m$.
\item Algorithm~\ref{alg:proposed_solver}.
\item Algorithm~\ref{alg:dual_fast_gradient}. 
\end{enumerate}
\begin{table}[t]
	\begin{center}
	\caption{Iterations of Algorithms~\ref{alg:primal_dual_interior_point} (columns 1-2) and~\ref{alg:proposed_solver} (column 3).}
		\label{tab:iter}
		\begin{tabular}{|c|rr|r|}
		\hline
			 & Scenario 1 & Scenario 2 & Scenario 3 \\\hline
			\emph{Damped Newton Phase}& 6 & 16 & -- \\
			\emph{Pure Newton Phase} & 20 & 8 & 6 \\\hline
		\end{tabular}
	\end{center}
\end{table}
\begin{rem}
Scenarios 1 and 2 rely on two different initializations of the dual
variables to show the impacts that their initialization have on the behavior of
the solver. To the best of our knowledge, many off-the-shelf interior-point
solvers do not allow the user to access the dual variables for their initialization (e.g., MATLAB's \texttt{quadprog}). 
A default choice is the one proposed in Scenario 1 in which the dual variables are initialized to $\mathbf 1$. 
The warm-starting strategy proposed in Scenario~2 considers $\lambda_0$ close to
zero, which is the optimal value of the multipliers associated with the unconstrained LS problem. 
\end{rem}

\begin{figure}[t]
	\begin{center} 
		\includegraphics[width=1\columnwidth]{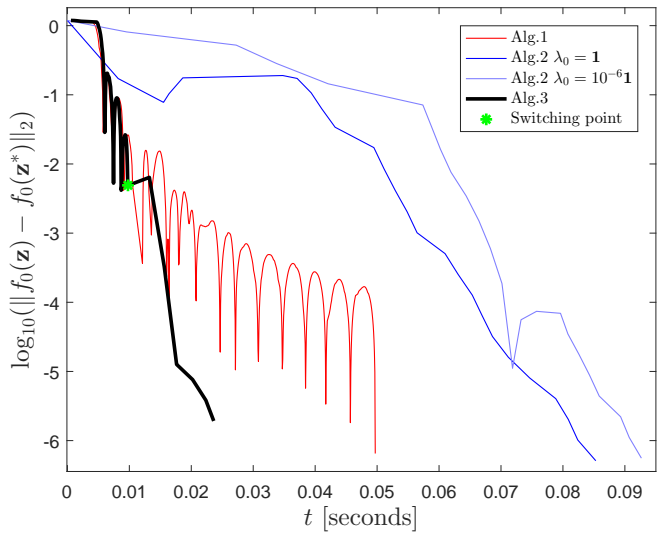}    
		\caption{Time required for online optimization using Algorithms~\ref{alg:dual_fast_gradient},~\ref{alg:primal_dual_interior_point}, and~\ref{alg:proposed_solver} for the 2-dimensional planar system.} 
		\label{fig:time_efficiency}
	\end{center} 
\end{figure}

Table~\ref{tab:iter} compares Algorithm~\ref{alg:proposed_solver} (Scenario 3) with Algorithm~\ref{alg:primal_dual_interior_point} (Scenarios 1 and 2) in terms of Newton iterations. Algorithm~\ref{alg:primal_dual_interior_point} requires 26 iterations in Scenario 1, while it requires 24 iterations in Scenario 2. Notice, however, that, while the total number of iterations in Scenario 2 is reduced, Algorithm~\ref{alg:primal_dual_interior_point} requires more damped-Newton-phase iterations. The main reason is that the initialization in Scenario 2 is farther from the optimal value, given that it is initialized with the pair $(\bfz_{\textrm{LS}},\lambda_0 \approx 0)$, which assumes no active constraints (while constraints are active at the optimum).
As Table~\ref{tab:iter} shows (Scenario 3), Algorithm~\ref{alg:proposed_solver} reduces the number of iterations to 6 and completely eliminates the damped Newton phase, thanks to the use of the DFG to initialize the interior-point iterates. This leads to significant improvements also from the computation point of view as Figure~\ref{fig:time_efficiency} depicts. In particular, Figure \ref{fig:time_efficiency} shows the computation time required to solve Problem~\eqref{eq:2_QP_MPC} online using the proposed algorithm. In Figure~\ref{fig:time_efficiency}, $f_0(\bfz^*)$ is the optimal cost corresponding to $\bfz^*$, computed by MATLAB's 
\texttt{quadprog}. PDIP terminates when the duality gap satisfies $\mu\le\varepsilon=10^{-6}$. The blue solid line shows the computation time of Algorithm~\ref{alg:primal_dual_interior_point} in Scenario 1. The light-blue solid line is the computation time of Algorithm~\ref{alg:primal_dual_interior_point} in Scenario 2. The red solid line shows the
computation time of Algorithm~\ref{alg:dual_fast_gradient} (Scenario 4).
The black bold line shows our proposed solver. The green star indicates the switching point when the 
solution returned by DFG satisfies condition~\eqref{eq:4_condition_modified} ($\approx 10^{-3}$ ). 
Note that, Scenario 4 requires 591 iterations of Algorithm~\ref{alg:dual_fast_gradient} to reach a high accuracy ($\approx 10^{-6}$), while Scenario 3 only requires 65 DFG iterates to reach the medium accuracy needed and enter the pure Newton phase in Algorithm~\ref{alg:proposed_solver}. Furthermore, note that the improvements in terms of computation time compared to Scenarios 1 and 2 are related to the following facts: (i) our algorithm does not require backtracking line search (computationally costly), and (ii) the damped Newton phase is replaced by DFG iterates. 
\begin{rem}
Figure \ref{fig:time_efficiency} monitors $\|f_0(\bfz_k)-f_0(\bfz^*)\|_2$ (in logarithmic scale). 
When Algorithm~\ref{alg:primal_dual_interior_point} is initialized with the optimal solution
$\bfz_0:=\bfz_{\textmd{LS}}$ (Scenarios 1 and 2), $f_0(\bfz_0)<f_0(\bfz^*)$,
given that it does not account for the presence of (active) constraints at the
optimum.
This leads to the nonmonotonic behavior at 0.02 sec (blue line) and 0.08
sec (light-blue line) in Figure~\ref{fig:time_efficiency}, when PDIP enters
the feasible region (and in pure Newton phase).
\end{rem}

\subsection{Cessna Citation Aircraft \cite{maciejowski2002predictive}}

The aircraft is flying at an altitude of 5000 km and a speed of 128.2 m/sec. The
system has 4 states, 1 input, 2 outputs and is subject to the following constrains: elevator angle $\pm 15^{\circ}$, elevator rate $\pm 30^{\circ}/$sec and pitch angle $\pm 30^{\circ}$. To discretize the system, we use a sampling time of $T=0.25$ sec. Moreover, $Q$ is identity matrix, $R=10$ and $N=10$. For this problem, $m_d=1.1394\times 10^{-4}$ and $L_d=5\times 10^{-7}$. Thus, $\eta_d = 2.6\times 10^{-2}$. The backtracking line search parameters for Algorithm \ref{alg:primal_dual_interior_point} are chosen as $\alpha=\frac{1}{3}$, $\beta=0.5$.

We test the proposed solver at 30 different initial states sampled uniformly in the feasible region $[0~60\text{km}]$ of altitude (other states fixed). For each state, the solver is run for 11 times and the median of the computational time is computed. In order to examine how much noise is in the measurements of computational time, the standard deviation is computed for each initial condition.
\begin{table}[]
	\centering
	\caption{Computational Time (sec) of Algorithms~\ref{alg:dual_fast_gradient},~\ref{alg:primal_dual_interior_point}, and~\ref{alg:proposed_solver}}
	\label{tab:aircraft}
	\begin{tabular}{|c|rr|r|}
		\hline
		& Alg.~\ref{alg:dual_fast_gradient} & Alg.~\ref{alg:primal_dual_interior_point}   & Alg.~\ref{alg:proposed_solver}      \\ \hline
		\emph{Best}    & 0.0063    & 0.0141 & 0.0024   \\ 
		\emph{Worst}   & --    & 0.0398 & 0.0214       \\ 
		\emph{Average} & --    & 0.0264 & 0.0120       \\ \hline
		 \begin{tabular}{@{}c@{}}\emph{Average} \\ \emph{Deviation}\end{tabular}& $1.4\times10^{-3}$ & $3.7\times10^{-4}$ & $3.6\times10^{-4}$ \\ \hline
	\end{tabular}
\end{table}

Table \ref{tab:aircraft} compares the three algorithms in terms of computational time. 
Each solver is supposed to reach the desired accuracy on duality gap $\mu\le\varepsilon=10^{-6}$. 
In the best case scenario, Algorithm \ref{alg:proposed_solver} (15 iterations
in DFG phase and 12 iterations in the pure Newton phase) only consumes $26.95\%$ of the
computational time of Algorithm \ref{alg:primal_dual_interior_point} 
(22 Newton iterations, 10 in the damped Newton phase), and $60.32\%$ of Algorithm \ref{alg:dual_fast_gradient} 
(280 iterations). In the worst case, Algorithm~\ref{alg:proposed_solver} (352
iterations in DFG phase and 11 iterations in the pure Newton phase) saves 46\%
of the computational time consumed by Algorithm \ref{alg:primal_dual_interior_point}
 (36 Newton iterations, 25 in damped Newton phase). Note that the
 execution of Algorithm \ref{alg:dual_fast_gradient} is terminated when it exceeded a given threshold (such as maximal number of iterations).
  On average, Algorithm \ref{alg:proposed_solver} saves more than 54\% of the
  computational time required by Algorithm \ref{alg:primal_dual_interior_point} thanks to elimination of 
  damped Newton phase.     
%
%
 \section{Conclusions and Future Work}
 \label{sec:conclusions}
This paper proposes an improved primal-dual interior point (PDIP) method for optimization problems that 
typically arise from model predictive control applications, that is, quadratic
programming problems with linear inequality constraints.
The proposed solver improves the convergence of state-of-the-art PDIP methods by replacing the damped Newton phase with a dual
 fast gradient method. This result is obtained by working in the dual space and modifying the condition 
 to enter the pure Newton phase. 
Finally, we showed the benefits of the proposed algorithm on a discrete unstable planar system and 
the Cessna Citation Aircraft system.

As part of the future work, we plan to test the technique by comparing with
 other solvers, such as, QPOASES, CPLEX, GUROBI.
%

\bibliographystyle{IEEEtran}
\bibliography{IEEEarxiv} 

\begin{thebibliography}{10}
\providecommand{\url}[1]{#1}
\csname url@rmstyle\endcsname
\providecommand{\newblock}{\relax}
\providecommand{\bibinfo}[2]{#2}
\providecommand\BIBentrySTDinterwordspacing{\spaceskip=0pt\relax}
\providecommand\BIBentryALTinterwordstretchfactor{4}
\providecommand\BIBentryALTinterwordspacing{\spaceskip=\fontdimen2\font plus
\BIBentryALTinterwordstretchfactor\fontdimen3\font minus
  \fontdimen4\font\relax}
\providecommand\BIBforeignlanguage[2]{{%
\expandafter\ifx\csname l@#1\endcsname\relax
\typeout{** WARNING: IEEEtran.bst: No hyphenation pattern has been}%
\typeout{** loaded for the language `#1'. Using the pattern for}%
\typeout{** the default language instead.}%
\else
\language=\csname l@#1\endcsname
\fi
#2}}

\bibitem{bemporad2002explicit}
A.~Bemporad, M.~Morari, V.~Dua, and E.~N. Pistikopoulos, ``The explicit linear
  quadratic regulator for constrained systems,'' \emph{Automatica}, vol.~38,
  no.~1, pp. 3--20, 2002.

\bibitem{nesterov1983method}
Y.~Nesterov, ``A method of solving a convex programming problem with
  convergence rate ${O}(1/k^2)$,'' in \emph{Soviet Mathematics Doklady},
  vol.~27, no.~2, 1983, pp. 372--376.

\bibitem{kogel2011fast}
M.~K{o}gel and R.~Findeisen, ``Fast predictive control of linear systems
  combining {Nesterov}'s gradient method and the method of multipliers,'' in
  \emph{Proc. of the 50th IEEE CDC}.\hskip 1em plus 0.5em minus 0.4em\relax
  IEEE, 2011, pp. 501--506.

\bibitem{nesterov2013introductory}
Y.~Nesterov, \emph{Introductory lectures on convex optimization: {A} basic
  course}.\hskip 1em plus 0.5em minus 0.4em\relax Springer Science \& Business
  Media, 2013, vol.~87.

\bibitem{Stathopoulos2016}
G.~Stathopoulos, H.~Shukla, A.~Szucs, Y.~Pu, and C.~N. Jones, ``Operator
  splitting methods in control,'' \emph{Foundations and Trends$^\circledR$ in
  Systems and Control}, 2016.

\bibitem{ricker1985ASMuse}
N.~L. Ricker, ``Use of quadratic programming for constrained internal model
  control,'' \emph{Industrial \& Engineering Chemistry Process Design and
  Development}, vol.~24, no.~4, pp. 925--936, 1985.

\bibitem{schmid1994ASMquadratic}
C.~Schmid and L.~T. Biegler, ``Quadratic programming methods for reduced
  hessian {SQP},'' \emph{Computers \& chemical engineering}, vol.~18, no.~9,
  pp. 817--832, 1994.

\bibitem{ferreau2006ASMonline}
H.~J. Ferreau, ``An online active set strategy for fast solution of parametric
  quadratic programs with applications to predictive engine control,''
  \emph{University of Heidelberg}, 2006.

\bibitem{rao1998application}
C.~V. Rao, S.~J. Wright, and J.~B. Rawlings, ``Application of interior-point
  methods to model predictive control,'' \emph{Journal of optimization theory
  and applications}, vol.~99, no.~3, pp. 723--757, 1998.

\bibitem{wang2010fast}
Y.~Wang and S.~Boyd, ``Fast model predictive control using online
  optimization,'' \emph{IEEE Transactions on {CST}}, vol.~18, no.~2, pp.
  267--278, 2010.

\bibitem{zeilinger2014realtimeRobustMPC}
M.~N. Zeilinger, D.~M. Raimondo, A.~Domahidi, M.~Morari, and C.~N. Jones, ``On
  real-time robust model predictive control,'' \emph{Automatica}, vol.~50,
  no.~3, pp. 683--694, 2014.

\bibitem{Borrelli2015predictive}
F.~Borrelli, A.~Bemporad, and M.~Morari, \emph{Predictive control for linear
  and hybrid systems}, 2015, available at
  \texttt{www.mpc.berkeley.edu/}\texttt{mpc-course-material}.

\bibitem{boyd2004convex}
S.~Boyd and L.~Vandenberghe, \emph{Convex optimization}.\hskip 1em plus 0.5em
  minus 0.4em\relax Cambridge university press, 2004.

\bibitem{maciejowski2002predictive}
J.~M. Maciejowski, \emph{Predictive control: with constraints}.\hskip 1em plus
  0.5em minus 0.4em\relax Pearson education, 2002.

\bibitem{nesterov2005smooth}
Y.~Nesterov, ``Smooth minimization of non-smooth functions,''
  \emph{Mathematical programming}, vol. 103, no.~1, pp. 127--152, 2005.

\bibitem{giselsson2014improved}
P.~Giselsson, ``Improved fast dual gradient methods for embedded model
  predictive control,'' \emph{IFAC Proceedings Volumes}, vol.~47, no.~3, pp.
  2303--2309, 2014.

\bibitem{necoara2014rate}
I.~Necoara and V.~Nedelcu, ``Rate analysis of inexact dual first-order methods
  application to dual decomposition,'' \emph{IEEE Transactions on Automatic
  Control}, vol.~59, no.~5, pp. 1232--1243, 2014.

\bibitem{necoara2015adaptive}
I.~Necoara, L.~Ferranti, and T.~Keviczky, ``An adaptive constraint tightening
  approach to linear model predictive control based on approximation algorithms
  for optimization,'' \emph{OCAM}, vol.~36, no.~5, pp. 648--666, 2015.

\bibitem{kerrigan2000soft}
E.~C. Kerrigan and J.~M. Maciejowski, ``Soft constraints and exact penalty
  functions in model predictive control,'' in \emph{International Conference
  (Control 2000), Cambridge}, 2000.

\bibitem{rubagotti2014stabilizing}
M.~Rubagotti, P.~Patrinos, and A.~Bemporad, ``Stabilizing linear model
  predictive control under inexact numerical optimization,'' \emph{IEEE {TAC}},
  vol.~59, no.~6, pp. 1660--1666, 2014.

\end{thebibliography}

\end{document}